\documentclass[12pt]{article}
\newtheorem{Th}{Theorem}\usepackage{amssymb}

\author{
Vladimir Blinovsky}
 
\date{Institute for Information Transmission Problems,\\
B. Karetnyi 19, Moscow, Russia,\\
vblinovs@yandex.ru}

\begin{document}
\title{Proof of Riemann hypothesis}
\maketitle

\hfill To my wife Luisi

\begin{center}
{\bf Abstract}\bigskip

 We prove the Riemann hypothesis
 \end{center}
 
We use the following definition of the zeta function: $z = \sigma + iT, \sigma \in (0,1)$  ([7]):
\begin{equation}
\label{s1}
\zeta (z) = \frac{1}{(1-2^{1-z})\Gamma (z)} \int_0^\infty \frac{x^{z-1}dx}{e^x+1} .
\end{equation}
The Riemann hypothesis states that all nontrivial zeroes of $\zeta (z)$  are concentrated on the line
$\sigma = 1/2$.
The Riemann hypothesis was a famous Hilbert problem (number eight of 23). It is also one
of the seven Clay Millennium Prize Problems. It was formulated in Riemanns 1859 Manuscript~\cite{1}.
Elegant, crisp, falsifiable, and far-reaching, this conjecture is the epitome of what a good conjecture
should be~\cite{2}.
It is the opinion of many mathematicians that the Riemann hypothesis is probably the most
important open problem in pure mathematics~\cite{3}.
In this paper we prove this hypothesis. We highlight the main steps of the proof.
\begin{itemize}\item
We use the integral representation~(\ref{s1})  of Riemann zeta function (RZ) in the open critical strip $\sigma\in (0,1)$.
\item We consider function $K(\sigma ,T)$ - the square of modulus of the integral in the definition of zeta function~(\ref{s1}) which has the
same set of zeroes on the open critical strip $(0,1)\times\mathbb{R}$ as RZ.
\item We make the transformations of the formula for $K(\sigma ,T)$.
\item We know that if zero of RZ in the critical strip has coordinates  $(\sigma ,T)$, then the point $(1-\sigma ,T)$
is also zero of RZ.
\item It is known~(\cite{b1}), that Riemann hypothesis is true for $T<3\cdot 10^{12}$. Next we assume that $T>3\cdot 10^{12}$. Function $K(\sigma,T)$ is symmetrical under inversion $T\to -T$, hence it is sufficient to consider the case of positive $T$.  We also note the following result
 from Ford (2002) [11].
\begin{Th}  If $z=\sigma + iT$  is a zero of the Riemann zeta function, then
\begin{equation}
\label{p1} \sigma < 1-\frac{1}{57.54(\ln (T))^{2/3}(\ln \ln (T))^{1/3}}.
\end{equation}
\end{Th}
Proof of this theorem uses Vinogradovs mean-value theorem [13]. Result in this last theorem have been improved in  several later papers, we will use  the follows~\cite{12}:
$$
\sigma <1-\frac{1}{6 \ln(T)}. 
$$
Hence there is no zeroes of RZ when
$\sigma  \in \Gamma=: \Biggl[0,  \frac{1}{6\ln (T)}\Biggr]\bigcup \Biggl[ 1-\frac{1}{6\ln (T)},1\Biggr]$
\item It is left to prove that second derivative of $K(\sigma ,T)$ is nonnegative   in the domain ${\cal D}=:([0,1]\setminus \Gamma,\ [3\cdot10^{12},\infty))$.
\item From above it follows that zeros of RZ in critical strip are possible at  $\sigma = 1/2$.
\item To prove that the second derivative of $K(\sigma,T)$ is positive in the  domain ${\cal D}$ we make further 
transformations of this second derivative and use the correlation
inequality.
\item This correlation inequality is FKG inequality and  we prove that function $f(h,T)$ in~(\ref{as1}) is monotone
decreasing.
\end{itemize}

\bigskip 

{\bf Proof.}
\bigskip

The square of the module of the integral in the definition of the zeta function is as follows:
$$
K(\sigma ,T) = \frac{1}{2}\int_0^\infty \int_0^\infty \frac{x^{z-1}y^{\bar{z}-1} +x^{\bar{z}-1}y^{z-1}}{(e^x+1)(e^y+1)}dx.
$$
Changing variable $y\to y^2/x$, we have
\begin{eqnarray*} && 
K(\sigma ,T) = 2\int_0^\infty\int_0^\infty \frac{y^{2\sigma -1}\cos\left(2T\ln\frac{x}{y}\right)}{x(e^x+1)(e^{y^2/x}+1)}dxdy.
\end{eqnarray*}

Changing variable $x\to ye^{ x/2}$ and using identity $\cos (x)=2\cos^2 (x/2)-1$ we have
\begin{eqnarray*}&&
K(\sigma ,T) =2\int_0^\infty \int_0^\infty \frac{y^{2\sigma -1}\cos (Tx)dxdy}{(e^{ye^{-x/2}}+1)(e^{ye^{x/2}}+1)}\\
&& =\frac{2}{T}\Biggl(2\int_0^\infty y^{2\sigma -1}\int_0^\infty \frac{\cos^2(x/2)dxdy}{(e^{ye^{-x/(2T)}}+1)(e^{ye^{x/(2T)}}+1)}\\
&&
-\int_0^\infty y^{2\sigma -1}\int_0^\infty \frac{dxdy}{(e^{ye^{-x/(2T)}}+1)(e^{ye^{x/(2T)}}+1)}\Biggr).
\end{eqnarray*}
Then
\begin{eqnarray} \label{dd1}
&& (K(\sigma ,T))^{\prime\prime}_{\sigma\sigma} =\frac{8}{T}\Biggl(2\int_0^\infty y^{2\sigma -1}(\ln (y))^2\int_0^\infty \frac{\cos^2 (x/2)dxdy}{(e^{ye^{-x/(2T)}}+1)(e^{ye^{x/(2T)}}+1)}\\
&&
- \int_0^\infty y^{2\sigma -1}(\ln (y))^2  \int_0^\infty \frac{dxdy}{(e^{ye^{-x/(2T)}}+1)(e^{ye^{x/(2T)}}+1)}\Biggr).\nonumber
\end{eqnarray} At first are going to prove the inequality
\begin{equation}
\label{d1}
(K(\sigma ,T))^{\prime\prime}_{\sigma\sigma} \geq 0,\ \sigma \in \left[\frac{1}{6\ln (T)}, 1-\frac{1}{6\ln (T)}\right],\ T>3\cdot10^{12}.
\end{equation}
RZ  is defined for $\sigma >1$  by the expansion  
$$
\zeta (z)=1+\sum_{n=1}^\infty \frac{1}{n^z}
$$
and has analytic continuation to the area $\sigma <1$ by functional equation~\cite{1}:
$$
\pi^{-z/2}\Gamma\left(\frac{z}{2}\right)\zeta (
z) =\pi^{-(1-z)/2} \Gamma\left(\frac{1-z}{2}\right)\zeta (1-z),\ \sigma <1.
$$
From the last formula it follows that with the root $\sigma_0  \in (0,1)$, function $K(\sigma ,T)$ has the root $1-\sigma_0$ and
we use the fact that the convex nonnegative function in ${\cal D}$ could not have more than one root in ${\cal D}$ and hence  this root can have coordinate only on the symmetry line $\sigma =1/2$. 

Using~(\ref{dd1})
we conclude that to prove the inequality~(\ref{d1}) it is sufficient  to prove the inequality
\begin{eqnarray*}&&
\int_0^\infty y^{2\sigma -1}(
\ln (y))^2 \int_0^\infty \frac{\cos^2 (\pi x/2)dxdy}{(e^{ye^{-\pi x/(2T)}}+1) (e^{ye^{\pi x/(2T)}}+1)}\\
&&\geq \frac{1}{2 }\int_0^\infty y^{2\sigma -1}(
\ln (y))^2 \int_0^\infty \frac{dxdy}{(e^{ye^{-\pi x/(2T)}}+1) (e^{ye^{\pi x/(2T)}}+1)}.
\end{eqnarray*}

The last inequality is equivalent to the inequality
\begin{equation}
\label{jj1}
\int_0^1 f(h,T)\cos^2 (\pi h/2)dh \geq\frac{1}{2}\int_0^1 f(h,T)dh,
\end{equation}
where
\begin{equation}
\label{as1}
f(h,T) =\sum_{i=0}^\infty \int_0^\infty y^{2\sigma -1}(\ln (y))^2 \bar{\lambda}(i,y,h,T)dy,
\end{equation}
\begin{eqnarray}&&
\bar{\lambda}(i,y,h,T) \\ \label{kk1}
&&= \frac{1}{(e^{ye^{-\pi i/T -\pi h/(2T)}}+1)(e^{ye^{\pi i/T +\pi h/(2T)}}+1)}+\frac{1}{(e^{ye^{-\pi (i+1)/T +\pi h/(2T)}}+1)(e^{ye^{\pi ( i+1)/T -\pi h/(2T)}}+1)}. \nonumber
\end{eqnarray}
Next we use  FKG inequality(see~\cite{3}: for two non-increasing    function $f,g: [0,1] \to \mathbf{R}_+$  the following correlation inequality is valid
\begin{equation}\label{ee1}
\int_0^1 f(h)g(h)dh\geq \int_0^1 f(h)dh \int_0^1 g(h)dh.
\end{equation}

Because function $g(h)=\cos^2 (\pi h/2)$ is monotone decreasing on $[0,1]$ and
$$
\int_0^1 \cos^2 (\pi h/2)=\frac{1}{2},
$$
due to correlation inequality~(\ref{ee1}) for the inequality~(\ref{jj1}) to be true it is suﬃcient to show that $f(h,T)$
decreasing with $h\in [0,1]$ or
$$
(f(h,T))^{(1)}_h \leq 0,\ h\in [0,1].
$$
We write 
$$
\bar{\lambda} (i,y,h,T) = \rho (i,y,h,T) + \rho (i,y,2-h,T),\ h\in [0,1],
$$
where
$$
 \rho (x,y,h,T) = \frac{1}{(e^{ye^{\pi x/T+\pi h/(2T)}}+1)(e^{ye^{-\pi x/T -\pi h/(2T)}}+1)}.
 $$
Thus to prove inequality~(\ref{jj1}) we need to prove inequality 
\begin{eqnarray}
\label{dsd1}&&
\left(\int_0^\infty y^{2\sigma-1}(\ln (y))^2\sum_{i=0}^\infty \bar{\lambda} (i,y,h,T) dy\right)^{(1)}_h \\&&= \int_0^\infty y^{2\sigma -1}(\ln (y))^2\sum_{i=0}^\infty ((\rho (i,y,h,T))^{(1)}_h - (\rho (i,y,2-h,T))^{(1)}_{2-h}) dy \leq 0. \nonumber
\end{eqnarray}
Rest of the article devoted to the proof of this inequality in the domain ${\cal D}$. 

 Define
$$
\rho (x,y,T)=\frac{1}{(e^{ye^{x\pi/T}}+1)(e^{ye^{-x\pi/T}}+1)}.
$$
We have 
$$
\rho (x+h/2,y,T)=\rho (x,y,h,T));\ (\rho (x,y,T))^{(1)}_{x;\ x=x+h/2}=2(\rho (x,y,h,T))^{(1)}_{h}.
$$
We use
Euler-Maclaurin  identity
$$
\sum_{i=0}^\infty g(z) = \int_0^\infty g(z)dz +\frac{1}{2} g(0) +\frac{1}{12}(g(z))^\prime_{z; z=0}+\frac{1}{6}\int_0^\infty B_3 (z-\lfloor x\rfloor ) (g(z))^{(3)}_{z} (z)dz,
$$
where
$$
B_3 (z)=z^3 -\frac{3}{2}z^2 +\frac{1}{2}z,\ |B_3 (z) |\leq  \frac{1}{2},\ z\in [0,1]; 
$$
 $g(x)\in C^\infty [0,\infty ]$ is smooth function and $g(\infty )=0$. We set $f(x,h,y,T):= (\rho (z,y,T))^{(1)}_{z:\ z=x+h/2} $ and $g(x,h,T):=\int_0^\infty y^{2\sigma-1}(\ln (y))^2(f(x,h,y,T)-f(x,2-h,y,T))dy$.  
 
Thus we have
\begin{eqnarray}&&\label{po9}
\Delta (y,T):= \sum_{i=0}^\infty (\rho^{(1)}_{z;\ z=i+h/2} (x,y,T)-\rho^{(1)}_{z;\ z=i +1-h/2}(z,y,T)) = \rho(2-h,y,T) -\rho(h,y,T) \\&&\nonumber \frac{1}{2}(\rho^{(1)}_{z; z=h/2} (z,y,T)-\rho^{(1)}_{z; z=1-h/2}(z,y,T))+\frac{1}{12}\nonumber  (\rho^{(2)}_{z; z=h/2} (z,y,T)-\rho^{(2)}_{z; z=1-h/2}(z,y,T))\nonumber \\&& +\frac{1}{6}\int_0^\infty B_3 (z-\lfloor z\rfloor )(\rho^{(4)}_{z; z=x+h/2} (z,y,T)-\rho^{(4)}_{z; z=x+1-h/2}(z,y,T))dx.\nonumber
\end{eqnarray}
Define
$$\Delta (T)=\int_0^\infty y^{2\sigma -1}(\ln (y))^2 \Delta (
y,T)dy.$$

Then 
\begin{eqnarray*}&&
\int_0^\infty y^{2\sigma -1}(\ln (y))^2\frac{dy}{(e^{ye^x}+1)(e^{ye^{-x}}+1)} = \int_0^\infty y^{2\sigma -1}(\ln (y))^2 \sum_{m ,n =1}^\infty (-1)^{m+n}e^{-y(n_1e^x+n_2e^{-x})}dy\\
&&\nonumber
=(\Gamma (2\sigma))^{(2)}\sum_{m ,n=1}^\infty \frac{(-1)^{m+n}}{(me^x+ne^{-x})^{2\sigma}} =\frac{1}{4}(\Gamma (2\sigma))^{(2)}\sum_{i=0}^\infty \frac{a_{2k}}{(2k)!} x^{2k},
\end{eqnarray*}
where
\begin{eqnarray*}&&a_{2k}=\frac{d^{2k}}{dx^{2k}} \sum_{m ,n=1}^\infty \frac{(-1)^{m+n}}{(me^x +n e^{-x})^{2\sigma}}_{x;\ x=0}\\&&=\frac{1}{4}\sum_{m,n=1}^\infty (-1)^{m +n}\sum_{\pi\in\Pi_{2k}}(-1)^{|\pi|}
\frac{(2\sigma)_{|\pi |}}{(m+n)^{2\sigma+|\pi|}}\prod_{B\in\pi}(m+(-1)^{|B|}n),
\end{eqnarray*}
where $\Pi_{2k}$ is the set of partitions of $2k-$ element set and $B\in \pi$ is element in the partition $\pi$  and $w(\pi),\ \pi\in\Pi_{2k}$ is the set of elements in $\pi$ with odd number of points. Note, that $|w(\pi)|=2u(\pi)$ is even. Then $|B|$ is size of element $B$ and $|\pi|$ is number of blocks in the partition $\pi$,\ $(x)_{m}=x(x+1)\ldots (x+m-1)\leq x (m+1)!$.

Note, that 
$\Gamma^{(2)}_{\sigma,\sigma}(2\sigma )=4\Gamma (2\sigma) ((\psi(2\sigma))^2+\psi^{(1)} (2\sigma))>0,\ \sigma\in (0,1)$, where $\psi (x)=(\ln\Gamma (x))^{(1)}$.
Define
\begin{eqnarray*}
&&
S(p,2\sigma ,x)=\sum_{m,n=1}^{\infty} (-1)^{m+n}\frac{(me^{x}-ne^{-x})^{p}}{(me^x+ne^{-x})^{2\sigma +p}},\ p\in {\bf N}\bigcup \{0\};\ \Omega(p,2\sigma):= S(p,2\sigma ,0).
\end{eqnarray*}
 We have
\begin{eqnarray*}&&
a_{2k}:=\sum_{m,n=1}^\infty (-1)^{m +n}\sum_{\pi\in\Pi_{2k}}(-1)^{|\pi|}
\frac{(2\sigma)_{|\pi |}}{(m+n)^{2\sigma+|\pi|}}\prod_{B\in\pi}(m+(-1)^{|B|}n)\\&&
= \sum_{q=0}^k\Omega(2q,2\sigma)\sum_{\pi\in\Pi_{2k}:\  \#\{ B\in\pi , : 2\not||B|\}=2q}(-1)^{|\pi|}(2\sigma )_{|\pi |}\\
&&
=\sum_{k_1=0}^k  {2k\choose 2k_1}\sum_{\pi_1\in\Pi_{2k_1,odd};\ \pi_2 \in\Pi_{2k-2k_1,even}} (-1)^{|\pi_2|}(2\sigma )_{|\pi_1|+|\pi_2|}\Omega(|\pi_1|,2\sigma)\\&&
=\sum_{k_1 =0}^{k}{2k\choose 2k_1}\sum_{q=0}^{k_1}\sum_{\pi_1\in\Pi_{2k_1 ,odd}, |\pi_1|=2q}\ \ \sum_{\pi_2\in\Pi_{2k-2k_1,even}}(-1)^{|\pi_2|}(x)_{2q+|\pi_2|}\Omega(2q,2\sigma).\end{eqnarray*}
Denote  Bernulli numbers $ B_{2j}$, then:
$$
\Omega(2q,2\sigma)=2^{1-2\sigma}\sum_{t=1}^\infty\frac{S_{2q}(t-1)}{t^{2\sigma +2q}}-2\sum_{t=1}^\infty \frac{S_{2q}(2t)-2^{2q}S_{2q}(t)}{(2t+1)^{2\sigma+2q}},
$$
where
$$
S_{2q}(n)=\sum_{k=1}^n k^{2q}=\frac{n^{2q+1}}{2q+1}+\frac{1}{2}n^{2q}+\sum_{j=1}^q {2q\choose 2j}\frac{B_{2j}}{2j}n^{2q-2j+1}.
$$
Define
$$G_2 (t,q) =2^{1-2\sigma}\sum_{t=1}^\infty\frac{S_{2q}(t-1)}{t^{2\sigma +2q}};\ G_1 (t,q)= 2\sum_{t=1}^\infty \frac{S_{2q}(2t)-2^{2q}S_{2q}(t)}{(2t+1)^{2\sigma+2q}}.$$

We have
\begin{eqnarray*}&&
S_{2q}(2t) =\frac{(2t)^{2q+1}}{2q+1}+\frac{1}{2}(2t)^{2q} + \sum_{j=1}^q {2q\choose 2j}\frac{B_{2j}}{2j}(2t)^{2q-2j+1};\\
&&
2^{2q}S_{2q}(t) =\frac{1}{2} \frac{(2t)^{2q+1}}{2q+1} +\frac{1}{2} (2t)^{2q}+\sum_{j=1}^q {2q\choose 2j}\frac{B_{2j}}{2j}(2t)^{2q-2j+1}2^{2j-1};\\
&&
S_{2q}(2t) -2^{2q}S_{2q}(t)=\frac{1}{2} \frac{(2t)^{2q+1}}{2q+1} +\sum_{j=1}^q {2q\choose 2j}\frac{B_{2j}}{2j}(2t)^{2q-2j+1}(1-2^{2j-1});\\
&&
G_1 (t,q):=\frac{S_{2q}(2t) -2^{2q}S_{2q}(t)}{(2t+1)^{2\sigma +2q}}=\frac{1}{(2t+1)^{2\sigma +2q}} \frac{(2t)^{2q+1}}{2(2q+1)} \\&& +\left(\frac{2t}{2t+1}\right)^{2\sigma+2q}\sum_{j=1}^q {2q\choose 2j}\frac{B_{2j}}{2j}(2t)^{-2\sigma-2j+1}(1-2^{2j-1});\\
&&S_{2q}(t-1)=S_{2q}(t)-t^{2q}=\frac{t^{2q+1}}{2p+1}-\frac{1}{2}t^{2q} +\sum_{j=1}^q {2q\choose 2j}\frac{B_{2j}}{2j}t^{2q-2j+1};\\
G_2 (t,q):= && \frac{S_{2q}(t-1)}{t^{2\sigma+2q}}=\frac{t^{1-2\sigma}}{2q+1}-\frac{1}{2}t^{-2\sigma}+\sum_{j=1}^q {2q\choose 2j}\frac{B_{2j}}{2j}t^{-2\sigma-2j+1}.
\end{eqnarray*}  
Hence
\begin{eqnarray*}&&
D(2q,2\sigma,t):=2^{1-2\sigma}G_2(t,q)-2G_1(t,q) =(2t)^{1-2\sigma}\frac{1}{2q+1}\left( 1-\left(\frac{2t}{2t+1}\right)^{2\sigma +2q}\right)-(2t)^{-2\sigma}\\
&& +\sum_{j=1}^q {2q\choose 2j}\frac{
B_{2j}}{2j}t^{-2\sigma-2j+2} + 2\left(\frac{2t}{2t+1}\right)^{2\sigma+2q}\sum_{j=1}^q {2q\choose 2j}\frac{B_{2j}}{2j}(2t)^{-2\sigma-2j+1}(2^{2j-1}-1),\\&&
\Omega(2q,2\sigma)=\sum_{t=1}^\infty D (2q,2\sigma ,t).
\end{eqnarray*}\
Because $q>0,$ we have 
\begin{eqnarray*}&&
|\Omega(2q,2\sigma )|=\Biggl|\sum_{t=1}^\infty 2^{1-2\sigma}G_2(t,q)-2G_1 (t,q)\Biggr|\leq 3\left(1+\sum_{j=1}^{q}{2q\choose 2j}\frac{|B_{2j}|}{2j}\right)\eta (2\sigma),\\
&&\leq 3(1+2^{2q}|B_{2q}|)\eta (2\sigma )<2^{2q+2}|B_{2q}|\zeta(2\sigma)\leq 2^{2q+2}(2q)!\eta (2\sigma),\  \sigma\in (0,1),\ \eta(x)=\sum_{i=1}^\infty \frac{(-1)^{n+1}}{n^x}, x >0.
\end{eqnarray*} 
Last inequality is true because Bernulli number $B_{2q}$ satisfies inequality $|B_{2(q-1)}|<|B_{2q}|<(2q)!,\ q>1$. 
Function $\eta (x)$ is Dirichlet  eta-function, $\eta (x)>0,\ x\in(0,2)$ is monotone increasing when $x>0$.

We need the following expansions 
\begin{eqnarray*}&&
(\sinh (x))^r =\frac{1}{2^r}\sum_{n=0}^\infty \frac{x^n}{n!} \sum_{i=0}^r (-1)^i {r\choose i}(2i-r)^n;\\&&(\sinh (x))^{2q} =2^{2n-2q}\sum_{n=0}^\infty \frac{x^{2n}}{(2n)!}\sum_{i=0}^{2q}{2q\choose i}(-1)^i(i-q)^{2n};\\&&
(\cosh (x)-1)^m =2^m(\sinh (x/2))^{2m}=2^{-m}\sum_{n=0}^\infty \frac{x^{2n}}{(2n)!} \sum_{k=0}^{2m} {2m\choose k}(-1)^k(k-m)^{2n};\\&&
\Delta (q,k,k_1)  :=\sum_{\pi\in\Pi_{2k_1}, |\pi|=2q, \ \forall A\in \pi, 2\not ||A|}\ \ \sum_{\xi\in\Pi_{2k-2k_1}, \forall B\in\xi, 2||A|} (-1)^{|\xi|}(x)_{2q+|\xi|}=(x)_{2q}\cdot O_{2k_1,2q}\\ && \times\sum_{m=0}^{k-k_1} E_{2k-2k_1,m}(-1)^m (x-2q)_m ;
\\&&
E_{2k-2k_1,m} =\frac{1}{m!}[x^{2k-2k_1}](\cosh (x)-1)^m =\frac{2^{-m+1}}{m! (2k-2k_1)!} \sum_{i=0}^{2m} {2m\choose i}(-1)^i (i-m)^{
2k-2k_1};\\&&
O_{2k_1,2q}=\frac{1}{(2q)!}[x^{2k_1}](\sinh (x))^{2q}=\frac{2^{2k_1-2q+1}}{(2k_1)!(2q)!}\sum_{i=0}^{2q}(-1)^i{2q\choose i}(i-q)^{2k_1}.
\end{eqnarray*}
Note, that $E_{2p,m} (O_{2p,2q})$ is the number of partitions of $[2p]$ into sets with $m$ even ($2q$ odd) number of elements. 
Thus
\begin{eqnarray*}&&
\frac{a_{2k}}{(2k)!}=\frac{1}{(2k)!}\sum_{k_1 =0}^k{2k\choose 2k_1}\sum_{q=0}^{k_1}  \Delta (q,k, k_1)  \Omega (2q,2\sigma)\\&&=\frac{1}{(2k)!^2}
 \sum_{m=0}^{k}\frac{1}{2^{m}m!}\sum_{j=0}^{2m} {2m\choose j}(-1)^j \ \sum_{q=0}^{k-m}\frac{1}{(2q)!}(x)_{2q+m} \Omega(2q,2\sigma )\sum_{i=0}^{2q}{2q\choose i}(-1)^i\\&&\times \sum_{k_1=q}^{k-m} {2k\choose 2k_1}^2 2^{2k_1}(j-m)^{2k-2k_1}(i-q)^{2k_1}
 \end{eqnarray*}
 and thus
 \begin{eqnarray*}
 &&\frac{|a_{2k}|}{(2k)!} \leq 2\sigma \frac{2^{8k}}{(2k)!^2}\sum_{m=0}^k \frac{(2q+m+1)! \Omega(2q,2\sigma )}{(2q)!m!}\cdot\left\{\begin{array}{ll} q^{2q}m^{2k-2q}, &m\geq q;\\ q^{2k-2m}m^{2m},& m<q \end{array}\right.\\&&
 \leq 8k \sigma\frac{2^{8k}}{(2k)!^2}\sum_{m=0}^k {2q+m\choose m}\Omega (2q,2\sigma)\cdot\left\{\begin{array}{ll} m^{2k}, &m\geq q;\\ q^{2k},& m<q\end{array}\right. < 2^{13k}\sigma\eta(2\sigma).
\end{eqnarray*}
Define
$$
\delta  (x\pi/T) =\int_{0}^\infty y^{2\sigma -1}(\ln(y))^2\rho (x,y,T)dy.
$$
We have
\begin{eqnarray}&&\label{rrr3}
\Delta (T) :=\delta (\pi(1-h/2)/T ))-\delta (\pi h/2T) +(\delta (\pi h/(2T)))^{(1)}_{h} +(\delta (\pi (1-h/2)/(T)))^{(1)}_h\\&&+\frac{2}{3}((\delta (\pi h/(2T)))^{(2)}_{h} -(\delta (\pi (1-h/2)/T))^{(2)}_h)\nonumber  \\ \nonumber
&&=\sum_{i=1}^\infty \frac{a_{2k}}{(2k)!}\left(\frac{\pi}{T}\right)^{2k} ((1-h/2)^{2k}-(h/2)^{2k})  -\sum_{k=1}^\infty \frac{a_{2k}}{(2k)!}(2k)\left(\frac{\pi}{T}\right)^{2k} ((1-h/2)^{2k-1}-(h/2)^{2k-1}) \\&&- \frac{2}{3}\sum_{k=2}^\infty \frac{a_{2k}}{(2k)!}(2k)(2k-1)\left( \frac{\pi}{T}\right)^{2k}(((1-h/2))^{2k-2}-(h/2)^{2k-2})\nonumber \\
&& \nonumber
=(1-h)\sum_{k=0}^\infty \frac{b_{2k}(h)}{(2k)!}\left(\frac{\pi}{2T}\right)^{2k}=-\frac{1}{24}(1-h)a_{4}(((1-h/2)h/2+2/3)\left(\frac{\pi}{T}\right)^4 -(1-h)\sum_{k=3}^\infty \frac{b_{2k}(h)}{(2k)!}\left(\frac{\pi}{T}\right)^{2k}
\end{eqnarray}
where
\begin{eqnarray*}&&
b_{2k}(h)= a_{2k}\left(\sum_{j=0}^{2k-1}(1-h/2)^j (h/2)^{2k-1-j} -\sum_{j=0}^{2k-2}(1-h/2)^j (h/2)^{2k-2-j} -\frac{2}{3}\sum_{j=0}^{2k-3}(1-h/2)^{j}(h/2)^{2k-3-j} \right) \\
&& =-a_{2k} \left(\sum_{j=0}^{2k-3} (1-h/2)^{j}(h/2)^{2k-3-j} ((1-h/2)h/2 +2/3)\right).
\end{eqnarray*} 
From above equation it follows that 
\begin{eqnarray*}&&
b_{2}(h)=0,\ b_4(h) =-a_{4}((1-h/2)h/2+2/3).
\end{eqnarray*}
Next we prove that when $\sigma\in (0,1)$, then
$$
a_{4} =(2\sigma )_{4}\Omega(4,2\sigma )-8\sigma(2\sigma +1)(3\sigma +2)\Omega (2,2\sigma )+4\sigma(3\sigma +1)\Omega(0,2\sigma )>0.
$$
Using inequalities
\begin{eqnarray*}&&
c_1 (a):=\min_{x\geq 2} x\left(1-\left(\frac{x}{x+1}\right)^a\right)-1 = 1-2\left(\frac{2}{3}\right)^{2\sigma};\\
&&c_2 (a):= \max_{x\geq 2}\frac{x}{3}\left(1-\left(\frac{x}{x+1}\right)^{a}\right)-1 =a-1;\\
&& c_3 (a):=\min_{x\geq 2}\frac{x}{5}\left(1-(\frac{x}{x+1})^{a}\right)-1 = -\frac{3}{5}-\frac{2}{5}\left(\frac{2}{3}\right)^a ,\ a\geq 1
\end{eqnarray*} 
 we have
\begin{eqnarray*}&&
\Omega (0,2\sigma ) = \sum_{t=1}^\infty \Biggl((2t)^{1-2\sigma}\left( 1-\left(\frac{2t}{2t+1}\right)^{2\sigma }\Biggr)-(2t)^{-2\sigma}\right)\\&& \geq 2^{-2\sigma} c_1 (2\sigma )\sum_{t=1}^\infty          t^{-2\sigma}\geq \left(1-2\left(\frac{2}{3}\right)^{2\sigma}\right)\eta (2\sigma)\geq  -\eta (2\sigma);\\
&&  \Omega(2,2\sigma)= \sum_{t=1}^\infty \Biggl((2t)^{1-2\sigma}\frac{1}{3}\left( 1-\left(\frac{2t}{2t+1}\right)^{2\sigma +2}\Biggr)-(2t)^{-2\sigma}\right)+t^{-2\sigma} + 2\left(\frac{2t}{2t+1}\right)^{2\sigma+2}(2t)^{-2\sigma-1}\Biggr)\\&&
\leq (c_2 (2\sigma+2)2^{-2\sigma} +1)\eta (2\sigma)+ \frac{1}{2}\eta (2\sigma)\leq \left({(2\sigma +1)2^{-2\sigma}}+3/2\right)\eta (2\sigma)\leq (2\sigma+5/2)\eta (2\sigma)\\&&
\Omega(4,2\sigma) =\sum_{t=1}^\infty \Biggl((2t)^{1-2\sigma}\frac{1}{5}\left( 1-\left(\frac{2t}{2t+1}\right)^{2\sigma +4}\right)-(2t)^{-2\sigma}\\
&& +\sum_{j=1}^2 {4\choose 2j}\frac{B_{2j}}{2j}t^{-2\sigma-2j+2} + 2\left(\frac{2t}{2t+1}\right)^{2\sigma+4}\sum_{j=1}^2 {4\choose 2j}\frac{B_{2j}}{2j}(2t)^{-2\sigma-2j+1}(2^{2j-1}-1)\Biggr)\\
&& = c_3 (2\sigma +4)2^{-2\sigma}\eta (2\sigma)+\sum_{t=1}^\infty \Biggl(-(2t)^{-2\sigma}\\ && + 6 t^{-2\sigma}+\frac{15}{4} t^{-2\sigma-2}+2\left(\frac{2t}{2t+1}\right)^{2\sigma+4}\left( 3\cdot 2^{-2\sigma} t^{-2\sigma -1}+\frac{105}{32} t^{-2\sigma -3}\right)\Biggr)\\
&& \geq -2^{-2\sigma}\eta(2\sigma )-2^{-2\sigma}\eta(2\sigma)+6\eta (2\sigma) +\frac{15}{4}\eta (2\sigma +2) +3 \cdot 2^{-2\sigma}\eta (2\sigma +1)+\frac{105}{32}\eta (2\sigma +3)
\\
&&
\geq -2^{-2\sigma+1}\eta(2\sigma )+6\eta (2\sigma) +\frac{15}{4}\eta (2\sigma)(1-2^{-2\sigma-2}) +3 \cdot 2^{-2\sigma}\eta (2\sigma )(1-2^{-2\sigma -1})+\frac{105}{32}\eta (2\sigma).\end{eqnarray*}
Here  we several times use inequality use inequality $1-2^{-x}<\eta (x)<1$. 
Define $b=2^{-2\sigma}\in (2^{-4},1)$, then we have 
\begin{eqnarray*}&&
\Omega (4,2\sigma) \geq \Biggl(-\frac{3}{2}b^2 -b \left(2+\frac{15}{16}+\frac{105}{256}-3\right)+6+\frac{15}{16}+\frac{105}{32}\Biggr)\geq 7\eta (2\sigma) .
\end{eqnarray*}
We have
\begin{eqnarray*}&&
a_4 \geq (2\sigma )_{4}\Omega(4,2\sigma )-8\sigma(2\sigma +1)(3\sigma +2)\Omega (2,2\sigma )+4\sigma(3\sigma +1)\Omega(0,2\sigma )\\&&>
2\sigma\Biggl( (2\sigma+1)(2\sigma+2)(2\sigma +3)\Biggl(-\frac{3}{2}b^2 -b \left(2+\frac{15}{16}+\frac{105}{256}-3\right)+6+\frac{15}{16}+\frac{105}{32}\Biggr)\\&& -4(2\sigma+1)(3\sigma +2)  (2\sigma +1)b+3/2)-2(3\sigma+1)\Biggr)\eta(2\sigma).
\end{eqnarray*}
Minimum in the rhs of the last chain of inequalities in $b$ achieved at $b=1$. Simplified expression in the rhs of last chain of inequalities we obtain expression ($x=2\sigma$, expression in brackets achieve  its minimum at $x=0$): 
\begin{eqnarray*}&&
a_4 > x(x^3+17x^2+45x+30)\eta(x)\geq 60\sigma\eta(2\sigma).\end{eqnarray*}
From other hand
\begin{eqnarray}&& \label{l1}
\frac{|b_{2k}(h)|}{(2k)!} \leq \frac{|a_{2k}|}{(2k)!}(1/4+2/3) \sum_{j=0}^{2j-3} 2^{-(2k-j-3)}\leq \frac{|a_{2k}|}{(2k)!}\frac{5}{12} \sum_{i=0}^\infty 2^{-i} \leq 2^{14k}\sigma\eta (2\sigma ).
\end{eqnarray}
Using Lagrange Theorem, we obtain the relation
\begin{eqnarray}&&\label{l2}
\Delta_2 (x,y): =\rho^{(4)}_{z;z=x+h/2} (z,y,T)-\rho^{(4)}_{z;\ z=x+1-h/2}(z,y,T)\\&&=-\frac{(1-h)\pi}{T} \rho^{(5)}(z,y,T)_{z;z=\xi};\ \hbox{for some}\ \xi\in [x+h/2,x+1-h/2].\nonumber
\end{eqnarray}
Defline
$$
\mu (x,y,T)= \frac{1}{e^{ye^x}+1}.$$
We need the upper bounds for the absolute values of derivatives $(\mu (x,y,0,T))^{(j)},\ j\in [5]$ and at first we write the explicit formulas for these derivatives: 

\begin{eqnarray*}
&&
(\mu (x,y,T))^{(1)}_x =-ye^{x}e^{ye^x}\frac{1}{(e^{ye^x}+1)^{2}};\\
&&
(\mu (x,y,T))^{(2)}_x=\frac{y(ye^x e^{ye^x}-ye^x-e^{ye^x}-1)e^x e^{ye^x}}{(e^{ye^x}+1)^3};\\&&
(\mu (x,y,T))^{(3)}_x =-ye^{x+ye^x}\frac{y^2e^{2x}e^{2ye^x}-4y^2 e^{2x}e^{ye^x}+y^2e^{2x}-3ye^xe^{2ye^x}+3ye^x +e^{2ye^x}+2e^{ye^x}+1}{(e^{ye^x}+1)^4};\\&&
(\mu (x,y,T))^{(4)}_x =4y^2e^{2x+2ye^x}\frac{y^2e^{2x}e^{2ye^x}-4y^2 e^{2x}e^{ye^x}+y^2e^{2x}-3ye^xe^{2ye^x}+3ye^x +e^{2ye^x}+2e^{ye^x}+1}{(e^{ye^x}+1)^5}\\
&& -ye^{x+ye^x}(1+ye^x)\frac{y^2e^{2x}e^{2ye^x}-4y^2 e^{2x}e^{ye^x}+y^2e^{2x}-3ye^xe^{2ye^x}+3ye^x +e^{2ye^x}+2e^{ye^x}+1}{(e^{ye^x}+1)^4}\\&&
-ye^{x+ye^x}\\&& \Biggl(2y^{2} e^{2x}e^{2ye^x} +2y^3 e^{3x}e^{2ye^x} -8y^2e^{2x}e^{ye^x}-4y^3 e^{3x}e^{ye^x}+2y^2 e^{2x}-3ye^x e^{2ye^x}-6y^2e^{2x}e^{2ye^x}\\&&
 +3ye^{x} +2ye^x e^{2ye^x}+2ye^xe^{ye^x}\Biggr)\frac{1}{(e^{ye^x}+1)^{4}}\\&&
(\mu (x,y,T))^{(5)}_x \\&&=8y^2e^{2x+2ye^x}(1+ye^x)\frac{y^2e^{2x}e^{2ye^x}-4y^2 e^{2x}e^{ye^x}+y^2e^{2x}-3ye^xe^{2ye^x}+3ye^x +e^{2ye^x}+2e^{ye^x}+1}{(e^{ye^x}+1)^5}\\&&
-20y^3 e^{3x+3ye^x}\frac{y^2e^{2x}e^{2ye^x}-4y^2 e^{2x}e^{ye^x}+y^2e^{2x}-3ye^x e^{2ye^x}+3ye^x +e^{2ye^x}+2e^{ye^x}+1}{(e^{ye^x}+1)^6}+
4y^2 e^{2x+2ye^x}\\&&\times \Biggl(2y^2e^{2x}e^{2ye^x} + 2y^3e^{3x}e^{2ye^x}-8y^2 e^{2x}e^{ye^x} - 4y^3 e^{3x}e^{ye^x} +2y^2 e^{2x}-3ye^x e^{2ye^x}-6y^2 e^{2x}e^{2ye^x}\\&&+3ye^x +2ye^x e^{2ye^x}+2ye^x e^{ye^x}\Biggr)\frac{1}{(e^{ye^x}+1)^{5}}\\&&
-(ye^{x+ye^x}(1+ye^x)^2 +y^2e^{2x+ye^x})\frac{y^2e^{2x}e^{2ye^x}-4y^2 e^{2x}e^{ye^x}+y^2e^{2x}-3ye^xe^{2ye^x}+3ye^x +e^{2ye^x}+2e^{ye^x}+1}{(e^{ye^x}+1)^4}\\&&
+4y^2 e^{2x+2ye^x}(1+ye^x)\frac{y^2e^{2x}e^{2ye^x}-4y^2 e^{2x}e^{ye^x}+y^2e^{2x}-3ye^xe^{2ye^x}+3ye^x +e^{2ye^x}+2e^{ye^x}+1}{(e^{ye^x}+1)^5}\\
&&
-2ye^{x+ye^x}(1+ye^x)\\&&\times\Biggl(2y^{2} e^{2x}e^{2ye^x} +2y^3 e^{3x}e^{2ye^x} -8y^2e^{2x}e^{ye^x}-4y^3 e^{3x}e^{ye^x}+2y^2 e^{2x}-3ye^x e^{2ye^x}-6y^2e^{2x}e^{2ye^x} \\&&+3ye^{x} +2ye^x e^{2ye^x}+2ye^xe^{ye^x}\Biggr)\frac{1}{(e^{ye^x}+1)^{4}}\\&&
+4y^2 e^{2x+2ye^x}\\&&\times\Biggl(2y^{2} e^{2x}e^{2ye^x} +2y^3 e^{3x}e^{2ye^x} -8y^2e^{2x}e^{ye^x}-4y^3 e^{3x}e^{ye^x}+2y^2 e^{2x}-3ye^x e^{2ye^x}-6y^2e^{2x}e^{2ye^x} +3ye^{x} \\&&+2ye^x e^{2ye^x}+2ye^xe^{ye^x}\Biggr)\frac{1}{(e^{ye^x}+1)^{5}}
-ye^{x+ye^x}
\Biggl(4y^{2} e^{2x}e^{2ye^x} + 4y^{3} e^{3x}e^{2ye^x} -16y^2e^{2x}e^{ye^x}\\&& -8y^3e^{3x}e^{ye^x}-12y^3 e^{3x}e^{ye^x}-4y^4 e^{4x}e^{ye^x}+4y^2 e^{2x}-3ye^x e^{2ye^x}  -6y^2e^{2x} e^{2ye^x}-12y^3e^{3x}e^{2ye^x}\\&& +3ye^x+2ye^xe^{2ye^x}+4y^2e^{2x}e^{2ye^x}
+2ye^xe^{ye^x}+2y^2 e^{2x}e^{ye^x}\Biggr)\frac{1}{(e^{ye^x}+1)^4}
\end{eqnarray*}  
Define
$$
\bar{\mu}_j (x,y,T))= (e^{ye^x}+1)|(\mu (x,y,T))^{(j)}_x|
$$
and $u=ye^x,$ then after some cumbersome calculations from previous expressions for derivatives f $\mu (x,y,0,T)$.  we obtain bounds (bounds are rather redundant):
\begin{eqnarray*}&&
\bar{\mu}_0 (x,y,T)\leq 1;\\
&&
\bar{\mu}_1 (x,y,T)\leq u ;\\
&&
\bar{\mu}_2 (x,y,T)\leq 2u(u+1);\\
&&
\bar{\mu}_3 (x,y,T)\leq 2u(3u^2 + 3e^u+2)\leq 6u(u^2+u+1);\\
&&
\bar{\mu}_4 (x,y,T))\leq u(24u^3 +24u^2+16u + (1+u)(6u^2+6u+4) +6u^3 +16u^2 +10u)\\&&
\leq 2u(23u^3+26u^2+18u+2)\leq 52u(u^3+u^2+u+1)\\
&&
\bar{\mu}_5 (x,y,T))\leq 400u(u^4+u^3+u^2+u+1).
\end{eqnarray*}
Define
\begin{eqnarray*}&&
C_1(x,y,T)=(\mu (x,y,T))^{(5)}_x\mu (-x,y,T);\ C_2(x,y,T)=(\mu (x,y,T))^{(4)}_x (\mu (-x,y,T))^{(1)}_x;\\&&C_3 (x,y,T)=(\mu(x,y,T))^{(3)}_x(\mu(-x,y,T))^{(2)}_x.
\end{eqnarray*}
Then
\begin{eqnarray*}&&
\left(\frac{1}{(e^{ye^x}+1)(e^{ye^{-x}}+1)}\right)^{5}_x =\sum_{i=0}^5 {5\choose i}(\mu (x,y,T))^{(i)}_x (\mu (-x,y,T))^{(5-i)}_x \\&&=C_1(x,y,T) - C_1 (-x,y,T) +5(C_2(x,y,T)-C_2 (-x,y,T))+10(C_3 (x,y,T)-C_3 (-x,y,T))  .
\end{eqnarray*}
We need the following transformation of $C_i (x,y,T)$:
\begin{eqnarray*} 
\bar{C}_i (x,y,T) := (e^{ye^x}+1)(e^{ye^{-x}}+1)C_i (x,y,T);\ i=1,2,3.
\end{eqnarray*}
Using Lagrange Theorem we obtain the inequality
\begin{eqnarray*}&&
|\bar{C}_1(x,y,T) - \bar{C}_1 (-x,y,T) +5(\bar{C}_2(x,y,T)-\bar{C}_2 (-x,y,T))+10(\bar{C}_3 (x,y,T)-\bar{C}_3 (-x,y,T))|\\&& \leq\frac{2\pi x}{T} \max_{\xi\in [x,-x]}|\bar{C}_1(\xi,y,T) +5\bar{C}_2(\xi,y,T)+10\bar{C}_3 (\xi,y,T)|.
\end{eqnarray*}
 Next we use formula for the derivative of product of two functions
$$
(f(x)g(x))^{(q)}_x =\sum_{i=0}^q {q\choose i}f^{(i)}_x (x)g^{(q-i)}_{x}(x),
$$
Wehave
\begin{eqnarray*}&&|C_1(x,y,T) - C_1 (-x,y,T) +5(C_2(x,y,T)-C_2 (-x,y,T))+10(C_3 (x,y,T)-C_3 (-x,y,T))|\\&&=\frac{1}{ (e^{ye^x}+1)(e^{ye^{-x}}+1)} |\bar{C}_1(x,y,T) - \bar{C}_1 (-x,y,T) +5(\bar{}C_2(x,y,T)-\bar{C}_2 (-x,y,T))\\&& +10(\bar{C}_3 (x,y,T)-\bar{C}_3 (-x,y,T))|. .
\end{eqnarray*}
Next we have
\begin{eqnarray*}&&
E_1: (x,y,T)=\max_{\xi\in [-x,x]}|\bar{C}_1 (
\xi,y,T) |=\max_{\xi\in [-x,x]}|(\mu (x,y,T))^{(5)}_x||\mu (-\xi,y,T)|\leq 400v(v^4+v^3+v^2+v+1);\\
&&
E_2 (x,y,T):=\max_{\xi\in [-x,x]}|C_2 (\xi,y,T) |=\max_{\xi\in [-x,x]}|(\mu (\xi,y,T))^{(4)}_x||(\mu (-\xi,y,T))^{(1)}_x|\leq 52y^2(v^3 +v^2+u+1);\\
&&
E_3 (x,y,T):=\max_{\xi\in [-x,x]}|C_3 (\xi,y,T)| =\max_{\xi\in [-x,x]}|(\mu (\xi,y,T))^{(3)}_x||(\mu (-\xi,y,T))^{(2)}_x|\leq 
24y^2(v^3+v^2+v+1),
\end{eqnarray*}
where $v=ye^{x/T},\ x\geq 0$.

Thus
\begin{eqnarray*}&&
\Delta (x,y):=\Biggl|\left(\frac{1}{(e^{ye^x}+1)(e^{ye^{-x}}+1)}\right)^{(5)}_x\Biggr| \leq \frac{2}{(e^{ye^x}+1)(e^{ye^{-x}}+1)} (E_1(x,y,T)+5E_2 (x,y,T)+10E_3 (x,y,T))\\&&
\leq e^{-v} ( 10^3 y^2(v^3+v^2+v+1) +400v(v^4+v^3+v^2+v+1)) <2\cdot 10^3 ve^{-v} \sum_{i=0}^4 v^i
\end{eqnarray*}

   Thus we have
 $$
 \Delta (x,y):=\max_{\xi\in [x+h/2,x+1-h/2]}   \Biggl|\Biggl(\frac{1}{(e^{ye^{x/T}}+1)(e^{ye^{-x/T}}+1)}\Biggr)^{(5)}\Biggr| \leq \frac{4x\pi^6}{T^6} \cdot 10^3 ve^{-v} \sum_{i=0}^4 v^i $$
 Next  we change variable $y\to ye^x$  to obtain the bound
 \begin{eqnarray*}&&
 \Delta_ 3(x,T):=\int_0^\infty y^{2\sigma-1}(\ln (y))^2 \Delta(x,y)dy \leq \frac{4x\pi^6}{T^6} \cdot 10^{3}\int_0^\infty y^{2\sigma -1}(\ln (y))^2  ve^{-v} \sum_{i=0}^4 v^i     dy         \\&& \leq \frac{4x\pi^6}{T^6} \cdot 10^{3}\int_0^\infty  v^{2\sigma } \left(\ln v  -\frac{x}{T}\right)^2e^{-v}  \sum_{i=0}^4 v^i             dv e^{-2x\sigma \pi/T}\\&&\leq\frac{8x\pi^6}{T^6}
  10^3\Biggl(\int_0^\infty v^{2\sigma } (\ln (v))^2e^{-v} \sum_{i=0}^4 v^i e^{-2\sigma \pi x/T}dz+\left(\frac{x\pi }{T} 
 \right)^2 \int_0^\infty z^{2\sigma } e^{-v} \sum_{i=0}^4 v^i e^{-2\sigma\pi x/T}\Biggr)dz\\&&\leq \frac{8x\pi^6}{T^6}
 10^3\Biggl( \frac{1}{4}\left(\sum_{i=1}^5\Gamma (2\sigma +i)\right)^{(2)}_{\sigma,\sigma } e^{-2\sigma \pi x/T} +\left(\frac{x\pi}{T}\right)^2 \sum_{i=1}^5 \Gamma (2\sigma +i)e^{-2\pi\sigma x/T}\Biggr). 
 \end{eqnarray*}
 Here in the second inequality we use inequality $(a+b)^2\leq 2(a^2+b^2)$.
 By Integration of the rhs of last chain of relations over $x\in (0,\infty )$ gives us the bound (here we use condition $\sigma >1/(6\ln( T))$:
 \begin{eqnarray*}&&
\Delta_4 (T) :=  \int_0^\infty \Delta_3 (x,T) dx\leq \frac{8\cdot 10^3 \pi^6}{T^6}
 \Biggl( \frac{1}{4}\sum_{i=0}^5(\Gamma (2\sigma +i))^{(2)}_{\sigma,\sigma } \int_0^\infty x e^{-2\sigma \pi x/T}dx \\&&+\frac{\pi^2}{T^{2}} \sum_{i=1}^5\Gamma (2\sigma +i)\int_0^\infty x^3 e^{-2\sigma  \pi x/T}dx\Biggr) \end{eqnarray*} 
Using relations
\begin{eqnarray*}&& 
 \psi (2\sigma +2) \leq\psi (z)=(\ln\Gamma (z))^{(1)} < \ln z<\ln 6,\ z\leq 6,\ (\psi )^{(1)}(z)\leq \frac{1}{z}+\frac{1}{2z^2}+\frac{1}{6z^3}<2,\\&& \ z>1;\ \Gamma (z)\leq  \Gamma (\lceil x\rceil )< \lceil x\rceil !,
 \Gamma (2\sigma+i)=\Gamma (2\sigma)(2\sigma)_i,
\end{eqnarray*}
we set inequalities
\begin{eqnarray*}&& 
\sum_{i=1}^5 \Gamma (2\sigma+i)<\Gamma(2\sigma)\sum_{i=1}^5 (2\sigma)_i=2\sigma\Gamma (2\sigma)\sum_{i=1}^4 (2\sigma +1)_{i}\\&&\leq 2\sigma\Gamma (2\sigma) \sum_{i=1}^4 \frac{ (2+i)!}{2}\leq 2\cdot 870\sigma\Gamma (2\sigma)\leq   1740\sigma;  \\&& \max_{x\in [0,7]}  \psi^{(1)}(x) <2, \   \max_{x\in [0,7]}  \psi (x) <\ln 6<2,\  \sum_{i=1}^5(\Gamma (2\sigma +i)(\psi^2 (2\sigma +i)+(\psi (2\sigma +i))^{(1)})\\&& \leq  6960\cdot 6\sigma<40000\sigma;\  
\end{eqnarray*}
and obtain the relations ($T>3\cdot10^{12}$):
\begin{eqnarray}&& \label{r1}
\Delta_4 (T) \leq \frac{8\pi^6}{ T^6 }10^3\Biggl(\frac{40000 T^2}{4\cdot4\sigma \pi^2}+\frac{\pi^2}{T^2}\frac{1740 \cdot 6 T^3}{(2 \pi)^3 \sigma^2}\Biggr)\\&&  \nonumber
 \leq \frac{2\pi^4}{ T^4\sigma }10^6\Biggl(\frac{1}{4}+\frac{17.4\cdot\pi\cdot 3\cdot6\ln (T)}{T}\Biggr)\\&& \nonumber
  \leq \frac{2\pi^4 \ln (T)}{ T^4 }10^6\Biggl(1+\frac{600\ln (T)}{T}\Biggr)  \leq \frac{\pi^4 }{ T^4 }10^{8}.
\end{eqnarray}
In the last chain of inequalities we use the inequality $T>3\cdot 10^{12}$.

Using formula~(\ref{l2})  and bounds~(\ref{l1}),~(\ref{r1}) we obtain bound:
\begin{eqnarray}&&\label{po99}
\Phi :=\int_0^\infty y^{2\sigma -1}(\ln (y))^2\sum_{i=0}^\infty \int_0^\infty\Biggl((\rho^{(1)}_{z;\ z=i+h/2} (x,y,T)-\rho^{(1)}_{z;\ z=i +1-h/2}(z,y,T))\Biggr)dydx \\&& \nonumber \leq \Delta (T) + \frac{1-h}{12T}\Delta_4 (T)\\ \nonumber
&&\leq -(1-h)\Biggl(\frac{1}{24}\left(\frac{\pi}{T}\right)^4 60\frac{2}{3}\sigma\eta (2\sigma) - \sum_{k=3}^\infty \left(\frac{2^{7}\pi}{T}\right)^{2k} \left(\frac{2}{3}+\frac{1}{4}\right)\sigma\eta (2\sigma)\\&& -  \frac{ \pi^4 }{ 12T^5 }10^{8}\Biggr)\Biggr)\\ \nonumber
&& \leq -(1-h) \frac{5\pi^4\sigma}{3T^4}\Biggl( \eta(2\sigma ) -\frac{3}{5}\frac{1}{1-\frac{\pi 2^{7}}{T}}\left(\frac{2^{7}\pi}{T}\right)^{2}\eta (2\sigma) -\frac{1}{20T}10^{10}\Biggr)\\&& \nonumber
 \leq -(1-h) \frac{5\pi^4}{3T^4}\Biggl( \eta(2\sigma )\left(1 -\frac{1}{10^{10}}\right) -.01\Biggr).
 \end{eqnarray}
 Because $\max_{x\in (0,2)}\eta (x)\geq .5$, we have the final bound
 $$
 \Phi   \leq -(1-h) \frac{\pi^4}{3T^4} . $$
 This prove the validness of the inequality~(\ref{d1}) and hence the convexity of $K(\sigma,T)$ over $\sigma\in [1/(6\ln (T)), 1-1/(6\ln(T))],\ T>3\cdot 10^{12}$ and from our considerations above   validness Riemann hypothesis follows.

 \end{document}